\documentclass[11pt]{article}
\usepackage{amsmath}
\usepackage{amssymb}
\usepackage{amsthm}
\usepackage[usenames]{color}
\usepackage{amscd}

\def\N{{\Bbb N}}
\def\je{{\bf j}}

\newtheorem{theorem}{Theorem}

\begin{document}

\title{\bf A note on the number of Abelian groups of a given order}
\author{L\'aszl\'o T\'oth \thanks{The author gratefully acknowledges support from the Austrian Science
Fund (FWF) under the project Nr. M1376-N18.}
\\ Institute of Mathematics, Department of Integrative
Biology \\ Universit\"at f\"ur Bodenkultur, Gregor Mendel-Stra{\ss}e
33, A-1180 Wien, Austria \\ and \\ Department of Mathematics,
University of P\'ecs \\ Ifj\'us\'ag u. 6, H-7624 P\'ecs, Hungary \\
E-mail: ltoth@gamma.ttk.pte.hu }
\date{}
\maketitle

\vskip-3mm

\centerline{{\it Mathematica Pannonica} {\bf 23} (2012), 157--160}

\vskip4mm

\centerline{\bf Dedicated to the memory of Professor Gyula I. Maurer (1927--2012)}

\begin{abstract} We point out an asymptotic formula for the
power moments of the function $a(n)$, representing the number of
non-isomorphic Abelian groups of order $n$. For the quadratic moment
this improves an earlier result due to L.~Zhang, M.~L\"u and
W.~Zhai.
\end{abstract}



\vskip3mm

Let $a(n)$ denote the number of non-isomorphic Abelian groups of
order $n$. The arithmetic function $a$ is multiplicative and for
every prime power $p^\nu$ ($\nu\ge 1$), $a(p^\nu)=P(\nu)$ is the
number of unrestricted partitions of $\nu$. Thus, for every prime
$p$, $a(p)=1$, $a(p^2)=2$, $a(p^3)=3$, $a(p^4)=5$, $a(p^5)=7$, etc.
Asymptotic properties of the function $a$ were investigated by
several authors. See, e.g., \cite[Ch.\ 14]{Iv1985}, \cite[Ch.\
7]{Kr1988} for historical surveys.

It is known that
\begin{equation*}
\sum_{n\le x} a(n) = A_1 x + A_2 x^{1/2} + A_3 x^{1/3} + R(x),
\end{equation*}
where $A_j:= \prod_{k=1, k\ne j}^{\infty} \zeta(k/j)$ ($j=1,2,3$),
$\zeta$ denoting the Riemann zeta function, and the best result for
the error term is $R(x)\ll x^{1/4+ \varepsilon}$ for every
$\varepsilon>0$, proved by O.~Robert and P.~Sargos \cite{RS2006}.
The asymptotic behavior of the sum $\sum_{n\le x} 1/a(n)$ was
investigated by W.~G.~Nowak \cite{No1991}.

An asymptotic formula for the quadratic moment of the function $a$,
i.e., for $\sum_{n\le x} (a(n))^2$ was given by L.~Zhang, M.~L\"u
and W.~Zhai \cite{Zh2008}.

In the present note we point out the following result for the $r$-th
power moment of the function $a$.

For $\je=(j_1,\ldots,j_t)\in \N^t$ with $1\le j_1\le \ldots \le j_t$
consider the generalized divisor function $d(\je;n):=
\sum_{d_1^{j_1}\cdots d_t^{j_t}=n} 1$ and let $\Delta(\je;x)$ stand
for the remainder term in the related asymptotic formula, i.e.,
\begin{equation*}
\sum_{n\le x} d(\je;n) = H(\je;x) + \Delta(\je;x),
\end{equation*}
where $H(\je;x)$ is the main term, cf. \cite[Ch.\ 6]{Kr1988}.
Furthermore, let $\Delta_r(x):=
\Delta((1,\underbrace{2,2,\ldots,2}_{2^r-1});x)$.

\begin{theorem} \label{theorem_1} Let $r\ge 2$ be a fixed integer. Assume that
$\Delta_r(x) \ll x^{\alpha_r} (\log x)^{\beta_r}$, with $1/3< \alpha_r < 1/2$. Then
\begin{equation*}
\sum_{n\le x} (a(n))^r = C_r x + x^{1/2} Q_{2^r-2}(\log x) + R_r(x),
\end{equation*}
where
\begin{equation*}
C_r: =\prod_p \left(1 + \sum_{\nu=2}^{\infty}
\frac{(P(\nu))^r-(P(\nu-1))^r}{p^\nu} \right),
\end{equation*}
$Q_{2^r-2}$ is a polynomial of degree $2^r-2$ and
$R_r(x)\ll x^{\alpha_r} (\log x)^{\beta_r}$ (is the same).
\end{theorem}

According to a recent result of E.~Kr\"atzel \cite{Kr2010},
$\Delta_2(x) \ll x^{45/127}(\log x)^5$, where $45/127\approx 0,3543
\in (1/3,1/2)$, hence the same is the remainder term for $\sum_{n\le
x} (a(n))^2$. This improves $R_2(x)\ll x^{96/245+\varepsilon}$ with
$96/245 \approx 0,3918$, obtained in \cite{Zh2008} by reducing the
error term to the Piltz divisor problem concerning $d_3(n)$.

If $r\ge 3$, then $\Delta_r(x)\ll x^{u_r+\varepsilon}$ for every
$\varepsilon>0$, where $u_r:=\frac{2^{r+1}-1}{2^{r+2}+1}\in
(1/3,1/2)$. See \cite[Th.\ 6.10]{Kr1988}. Therefore $R_r(x)\ll
x^{u_r+\varepsilon}$ holds as well.

Theorem \ref{theorem_1} is a direct consequence of the next more
general result, valid for a whole class of arithmetic functions. Let
$\Delta_{k,\ell}(x):=
\Delta((1,\underbrace{\ell,\ell,\ldots,\ell}_{k-1});x)$.

\begin{theorem} \label{theorem_2} Let $f$ be a complex valued multiplicative
arithmetic function. Assume that

i) $f(p)=f(p^2)=\cdots =f(p^{\ell-1})=1$, $f(p^{\ell})=k$ for every
prime $p$, where $\ell, k\ge 2$ are fixed integers,

ii) $f(p^\nu) \ll 2^{\nu/(\ell+1)}$ ($\nu \to \infty$) uniformly for
the primes $p$, i.e., there is a constant $C$ such that
$|f(p^\nu)|\le C\cdot 2^{\nu/(\ell+1)}$ for every prime $p$ and
every sufficiently large $\nu$.

Then
\begin{equation*}
\sum_{n=1}^{\infty} \frac{f(n)}{n^s} = \zeta(s)\zeta^{k-1}(\ell s)
V(s),
\end{equation*}
absolutely convergent for $\Re (s) > 1$, where the Dirichlet series
$V(s)$ is absolutely convergent for $\Re (s) > 1/(\ell+1)$.

Furthermore, suppose that $\Delta_{k,\ell}\ll x^{\alpha_{k,\ell}}
(\log x)^{\beta_{k,\ell}}$, with $1/(\ell+1) < \alpha_{k,\ell}<
1/\ell$. Then
\begin{equation*}
\sum_{n\le x} f(n)= C_f x + x^{1/\ell} P_{f,k-2}(\log x) + R_f(x),
\end{equation*}
where $P_{f,k-2}$ is a polynomial of degree $k-2$,
\begin{equation*}
C_f: =\prod_p \left(1 + \sum_{\nu=\ell}^{\infty}
\frac{f(p^\nu)-f(p^{\nu-1})}{p^\nu} \right),
\end{equation*}
and $R_f(x)\ll x^{\alpha_{k,\ell}} (\log x)^{\beta_{k,\ell}}$.
\end{theorem}

Note that for every $k,\ell \ge 2$, $\Delta_{k,\ell}(x)\ll
x^{u_{k,\ell}+\varepsilon}$, where $u_{k,\ell}: =
\frac{2k-1}{3+(2k-1)\ell}\in (1/(\ell+1),1/\ell)$. See \cite[Th.\
6.10]{Kr1988}. Therefore $R_f(x)\ll x^{u_{k,\ell}+\varepsilon}$ is
valid as well.

\begin{proof} This is a variation of the Theorem proved in
\cite{To2007}. Here the same proof works out, however the conditions
are somewhat relaxed. Let $\mu_{\ell}(n)=\mu(m)$ or $0$, according
as $n=m^{\ell}$ or not, where $\mu$ is the M\"obius function. Let
$V(s):=\sum_{n=1}^{\infty} v(n)/n^s$. We obtain the desired
Dirichlet series representation by taking $v=f*\mu *
\underbrace{\mu_{\ell}*\cdots \mu_{\ell}}_{k-1}$  in terms of the
Dirichlet convolution $*$.

Here $v$ is multiplicative and easy computations show that
$v(p^\nu)=0$ for any $1\le \nu \le \ell$. For $\nu \ge \ell +1$,
\begin{equation*}
v(p^\nu)= \sum_{j=0}^{k-1} (-1)^j {k-1 \choose j}
\left(f(p^{\nu-j\ell})-f(p^{\nu-j\ell-1})\right),
\end{equation*}
leading to the absolute convergence of $V(s)$ for $\Re (s)>1/(\ell +1)$.
Now the asymptotic formula follows from the representation
\begin{equation*}
f(n)= \sum_{ab=n} d((1,\underbrace{\ell,\ell,\ldots,\ell}_{k-1});a)v(b).
\end{equation*}
\end{proof}

Choosing $f(n)=(a(n))^r$, $k=2^r$ and $\ell=2$ we deduce Theorem
\ref{theorem_1}. Note that $P(\nu)<e^{\pi \sqrt{2\nu/3}}$ ($\nu \ge
1$), see e.g., \cite[p.\ 236]{Na1983}, thus condition ii) is
verified.

Theorem \ref{theorem_2} applies also for the $r$-th powers ($r\ge 2$
integer) of the exponential divisor function $\tau^{(e)}$ and the
function $\phi^{(e)}$, where $\phi^{(e)}$ is multiplicative and
$\phi^{(e)}(p^\nu)= \phi(\nu)$ for every prime power $p^\nu$ ($\nu
\ge 1$), $\phi$ denoting Euler's function. See \cite{Kr2010,
To2007}.

\end{document}